\documentclass[reqno,showkeys]{amsart}
\usepackage{amsfonts}
\usepackage{amssymb}
\usepackage{amsthm}
\usepackage{upref}
\makeatletter
\def\LaTeX{\leavevmode L\raise.42ex
    \hbox{\kern-.3em\size{\sf@size}{0pt}\selectfont A}\kern-.15em\TeX}
 \makeatother

\newcommand{\e}{\eqref}
\newcommand{\q}{\quad}
\newcommand{\ri}{\rightarrow}

\renewcommand{\S}{{\Bbb S}}
 
  \DeclareMathOperator*{\slim}{s-lim}

\numberwithin{equation}{section}

\newtheorem{lemma}{Lemma}[section]
\newtheorem{theorem}[lemma]{Theorem} 
\newtheorem{corollary}[lemma]{Corollary}
\newtheorem{proposition}[lemma]{Proposition}

\theoremstyle{definition}

\newtheorem{example}[lemma]{Example}

\newtheorem{assumption}[lemma]{Assumption}

\newtheorem{remark}[lemma]{Remark}

\def\qqq{\mathrel{\subset\mkern-15mu\lower.38ex\hbox{${\scriptscriptstyle\rightarrow}$}}}

\let\cal\mathcal

\let\Bbb\mathbb

\begin{document}

\title{Scattering by magnetic fields} 
\author{ D. R. Yafaev}
\address{ IRMAR, University of Rennes-I, Campus  
  Beaulieu, 35042 Rennes Cedex, FRANCE}
\email{yafaev@univ-rennes1.fr}
 
 \begin{abstract}
Consider the scattering amplitude $s(\omega,\omega^\prime;\lambda)$,
$\omega,\omega^\prime\in{\Bbb S}^{d-1}$, $\lambda > 0$,
corresponding to an arbitrary   short-range
magnetic field $B(x)$, $x\in{\Bbb R}^d$. This is a smooth function of $\omega$ and $\omega^\prime$ away from the diagonal $\omega=\omega^\prime$ but it may be singular on the diagonal. If $d=2$, then the singular part of the scattering amplitude
 (for example, in the transversal gauge) is a linear combination
of the Dirac function and of a singular denominator. Such structure is typical for long-range scattering. We refer to this phenomenon as to
the long-range Aharonov-Bohm effect. On the contrary, for $d=3$   
 scattering is essentially of short-range nature although, for example,
 the magnetic
 potential $A^{(tr)}(x)$  such that ${\rm curl}\,A^{(tr)}(x)=B(x)$ and
$ \langle A^{(tr)}(x),x\rangle=0$  decays  at infinity as $|x|^{-1}$ only. To be more precise,
we show that, up to  the diagonal Dirac function (times  an explicit function of $\omega$),   the scattering amplitude  has only  a weak singularity  
 in the forward direction $\omega = \omega^\prime$.
  Our approach relies on a construction in the dimension $d=3$  of a       short-range  magnetic  potential $A (x)$  corresponding to a given short-range  magnetic field $B(x)$.
 \end{abstract}

 \maketitle

\section{Introduction}

 {\bf 1.1.} 
In the original paper  \cite{AB} by Aharonov and Bohm (see also \cite{Hen, RUI}) the following mental experiment was discussed.  Consider a  thin straight solenoid
of  infinite length so that the magnetic field $B(x)$ is confined inside this solenoid and is zero outside of it.  Consider a beam of particles (electrons) coming from infinity following some direction. Suppose that its interaction with the magnetic field  inside the solenoid is blocked out by some shield, for example, by a strong repulsive electric field. Nevertheless the scattering amplitude turns out to be different from zero (the corresponding scattering matrix is not the identity operator). Therefore it can be expected that one may observe 
in experiments a non-trivial interference behind the solenoid between parts of the initial beam going around the solenoid from the left and right. Moreover, this
interference picture should
  depend on the magnetic flux $\Phi$ through  a cross-section of the solenoid. This contradicts of course the classical picture but is perfectly conformal with the principles of quantum mechanics.  Indeed, the Schr\"odinger equation is formulated in terms of a magnetic potential $A (x)$ defined by the equation
 \begin{equation}\label{eq:potfield}
  {\rm curl}\, A(x)=B(x). 
\end{equation}
  In view of translation invariance in the direction of 
the solenoid, the problem considered is two-dimensional. For definiteness, we suppose that the axis of the solenoid coincides with the $x_{3}$-axis, so that
$B(x)=(0,0, B(x))$, $x=(x_1,x_2)$, 
and (\ref{eq:potfield}) reduces to the equation
\begin{equation}\label{eq:potfield2}
  \partial A_{2}(x)/\partial x_{1}-  \partial A_{1}(x)/\partial x_{2}=B(x)
\end{equation} 
for components of the potential $A(x)=(A_{1}(x), A_{2}(x),0)$.
According to the Stokes theorem, the flux $\Phi$ is defined by one of the following two equalities
 \begin{equation}\label{eq:flux}
\Phi =\int_{{\Bbb R}^2} B(x) dx
=\lim_{R\rightarrow\infty}\int_{|x|=R}\langle A(x),dx\rangle 
\end{equation}
   where  $ \langle\cdot, \cdot\rangle $ is the scalar product in ${\Bbb R}^2$ or, more generally, in ${\Bbb R}^d$.
 Therefore  $A (x)$ is not zero   provided $\Phi\neq 0$. It follows that   scattering   is non-trivial     even though
the magnetic field   is   zero   in the region where particles can  penetrate.

Actually, to solve the problem explicitly, Aharonov and Bohm have simplified it in the following way. First, they   chose an infinitely  thin  solenoid. Second, instead of an inpenetrable shield described mathematically by the Dirichlet boundary condition they put a regular condition on the solenoid itself (at $x=0$). Thus, strictly speaking, a direct interaction of particles with a magnetic field was not   completely  excluded in the Aharonov-Bohm (A-B) model (this ``drawback" was remedied in \cite{RUI}, see also \cite{Af}). To be more precise,
in the   paper \cite{AB}   the Schr\"odinger operator $H$  
with magnetic potential 
\begin{equation}\label{eq:ABp}
A (x)=-\alpha (-x_2,x_1,0) |x|^{-2},
\quad\alpha= -(2\pi)^{-1}\Phi \in{\Bbb R}, 
\end{equation} 
 was considered (according to \e{eq:flux}, in this case the magnetic  field  
equals $- 2\pi \alpha \delta(x)$ where $\delta(x)$ is the
Dirac function).
 For such potentials, variables in  the Schr\"odinger equation
 can be separated  in polar coordinates $(r,\theta)$, and for
every   angular momentum $m=0, \pm 1, \pm 2,\ldots$ the radial equation 
\[
-u_m^{\prime\prime}+((m+\alpha)^2-1/4)r^{-2}u_m=\lambda u_m 
\]
($\lambda > 0$  is the energy)
can be solved in terms of the Bessel functions ${\cal
I}_{\nu}$, namely $u_m(r)=r^{1/2}{\cal
I}_{|m+\alpha|}(\lambda^{1/2}r)$. Using the asymptotics of these functions as $r\rightarrow\infty$, we see that
the scattering matrix (SM) $S $  for the  operator $H$ does not depend on $\lambda$ and has two eigenvalues
\begin{equation}\label{eq:EV}
s_{m}=e^{i\alpha\pi}\quad \mbox{for} \quad m\leq - \alpha  \qquad \mbox{and}
\qquad  s_{m}=e^{-i\alpha\pi} \quad \mbox{for} \quad m\geq - \alpha 
\end{equation} 
  with corresponding eigenfunctions $e^{im\theta}$. This implies of
course that the SM is non-trivial, that is, $S \neq I$ ($I$ is the identity operator) if $\alpha\not\in 2{\Bbb Z}$.
Since the magnetic field $B(x) =0$ for $x\neq 0$,
the fact that $S \neq I$, known now as the A-B effect,  appeared to be surprising,
at least from the point of view of classical physics. Anyway the A-B effect is a perfect test  for the validity  of quantum mechanics. Its experimental confirmation is discussed in \cite{PT}. Note that if we introduce the Planck constant $\hbar$, then $s_{m}(\hbar)=\exp(i\alpha \pi \hbar^{-2})$  for $m\leq - \alpha\hbar^{-2}$ and
$s_{m}(\hbar)=\exp(-i\alpha \pi \hbar^{-2})$ for $m\geq - \alpha \hbar^{-2}$. Therefore the SM $S(\hbar)$ has no limit as $\hbar\ri 0$ which is consistent with the absence of the A-B effect in classical physics.

Actually,  for the A-B 
potential, the SM is not only non-trivial, but its properties are typical for long-range scattering. Indeed, the scattering amplitude (kernel of the SM regarded as an integral operator) equals
\[
s(\theta,\theta^{\prime}) = (2\pi)^{-1}\sum_{m=-\infty}^\infty s_{m}e^{im(\theta - \theta^{\prime})}.
\]
A simple calculation (see \cite{RUI}) shows that, for eigenvalues (\ref{eq:EV}), this expression can be written as
\begin{equation}\label{eq:kern}
s(\theta,\theta^{\prime}) = \delta (\theta - \theta^{\prime}) \cos \pi\alpha+ i\pi^{-1}
  e^{- i[\alpha](\theta - \theta^{\prime})} \sin \pi\alpha \, P.V. (e^{i(\theta - \theta^{\prime})}-1)^{-1},
\end{equation} 
where $[\alpha]$ denotes the greatest integer less than or equal to $\alpha$. In particular,   the scattering amplitude contains a singular denominator (understood in the sense of the principal value) if the magnetic flux $\Phi\not\in 2\pi {\Bbb Z}$.
We use the term  ``long-range Aharonov-Bohm" effect for this phenomenon since  such singularity is absent for short-range potentials.

\bigskip

  {\bf 1.2.}
  Recall   that, under natural assumptions,
the scattering amplitude $s(\omega,\omega^{\prime};\lambda)$, where
$\omega,\omega^{\prime}\in{\Bbb S}^{d-1}$ if  $x\in{\Bbb R}^{d}$,  is a  smooth function away from the diagonal but can be very singular for $\omega=\omega^{\prime}$. This singularity is determined by the decay of a potential at infinity. For short-range   potentials    satisfying   the
condition 
\begin{equation} \label{eq:H1sr}
 |  A(x)|  \leq C  (1+|x|)^{- \rho }, \quad
 \rho >1,
\end{equation}
 the SM is the sum of
   the identity operator and of an integral operator with a weak diagonal singularity.
  Moreover, if $\rho\in (1,d)$, then the estimate
    \[
 s(\omega,\omega^{\prime}; \lambda)= O(| \omega - \omega^{\prime}|^{-d+\rho}),\quad
 \omega \neq\omega^{\prime}, \quad \omega - \omega^{\prime}\ri 0,
 \]
holds (see, e.g., \cite{Y5}). This   singularity is weaker than that of kernel of the singular integral operator (cf.  \e{eq:kern}). 
 Moreover, it   becomes weaker as long as a potential decays faster at infinity.
 
  For long-range potentials  the diagonal singularity
 of the scattering amplitude is stronger than in the short-range case, but the diagonal Dirac function disappears. Roughly speaking, for potentials $A(x)$  asymptotically homogeneous of degree $-\rho$, $\rho\in(0,1)$,  the   singular part  of the scattering amplitude is given by the formula
  \begin{equation} \label{eq:ARLR}
s_{0}(\omega,\omega^\prime;\lambda)=G(\omega,\omega-\omega^\prime;\lambda)
   \exp\Bigl(i\Xi (\omega,\omega-\omega^\prime;\lambda)\Bigr),
\end{equation}
   where $G$ and $\Xi$ are  asymptotically homogeneous functions of degrees $-(d-1)(1+\rho^{-1})/2$ and $1-\rho^{-1}$, respectively (see  \cite{Y1}). Thus,   
    the diagonal singularity of function (\ref{eq:ARLR}) is stronger than that of kernel of a   singular  integral operator. Nevertheless due to oscillations of the second factor in   (\ref{eq:ARLR}), the operator with such kernel is bounded in  $L_{2}({\Bbb S}^{d-1})$.
 
 Let us consider finally the intermediary case of    potentials with Coulomb decay ($\rho=1$) at infinity. Now  the results for electric and magnetic fields are qualitatively different. In the first case the answer is again given  (see, e.g., \cite{Her, Y}) by formula
     (\ref{eq:ARLR}), where $G$ is an asymptotically homogeneous function of degree $-d+1 $ and $\Xi$ has a logarithmic singularity at
  $\omega = \omega^{\prime}$. As shown in this paper, for arbitrary magnetic potentials satisfying the transversal condition
\begin{equation}\label{eq:TRG}
\langle  A(x),x\rangle =0, 
\end{equation}
the singularity of the scattering amplitude is described by a formular similar to
 (\ref{eq:kern}). Note that  we pay a special attention here
 to potentials corresponding to short-range magnetic fields, for which the answer depends crucially on  dimension of the space.      

 \bigskip

 {\bf 1.3.}
 Thus, we consider an arbitrary magnetic field satisfying the short-range condition
\begin{equation} \label{eq:B1}
 |  B(x)|  \leq C  (1+|x|)^{-r},\quad r > 2,
\end{equation}
($C$ denotes different positive constants whose precise values are of no importance)  and study properties of the corresponding SM. In particular, we find a diagonal singularity of the scattering amplitude. Our goal here is to reveal the difference between dimensions $2$ and $3$.
 
 Of course by the study of the SM for a given magnetic field one has to take into account that, from a theoretical point of view,  
the SM is  determined by a  magnetic potential $A(x)$ satisfying equation (\ref{eq:potfield}). In its turn, a solution of this equation is not unique, that is the gradient of an arbitrary function can be  added to $A(x)$ which leads to a gauge transformation of the Schr\"odinger operator. Although the SM are different in different gauges, they are connected by a simple formula, i.e., they are covariant with respect to gauge transformations. This allows us to speak about the SM corresponding to a given magnetic field.

As far as the A-B  effect  is concerned, the situation is similar in dimensions two and three. 
Consider, for example,   a toroidal solenoid
 ${\bf T}$ in the space ${\Bbb R}^3$.
  The magnetic field is  again concentrated inside the solenoid and is zero outside of it. Suppose again that this   solenoid is surrounded by a slightly bigger
toroidal solenoid which excludes a direct interaction of quantum particles with the magnetic field. By virtue of the Stokes theorem the corresponding magnetic potential is   non-zero and hence the SM is non-trivial provided the magnetic flux $\Phi_{s}$ through  a transverse cross-section of the solenoid is not zero.

On the contrary, it turns out that the long-range A-B  effect always occurs in dimension two (if the total magnetic flux $\Phi\not\in 2\pi {\Bbb Z}$), but under assumption \e{eq:B1} it cannot happen in dimension three. To be more precise, we show that in dimension
$d=2$, the diagonal  singularity of the scattering amplitude is described by a formula similar to (\ref{eq:kern}), i.e., it has a  structure typical for long-range scattering  (although the usual wave operators exist in this case). On the contrary, 
 in   dimension $d=3$ the structure of the SM is almost the same as for  scattering 
 by short-range potentials.
 We show also that condition 
(\ref{eq:B1}) is precise, that is the diagonal  singularity of the scattering amplitude is described by a formula generalizing (\ref{eq:kern}) to the case $d=3$, if
 $B(x)$ decays at infinity as a homogeneous function of degree  $-2$. Note that a condition similar to  (\ref{eq:B1}) distinguishes also short-range electric fields. 
 
A priori the difference between dimensions 
  $d=3$ and $d=2$ is not quite obvious. Indeed, in the case $d=3$ a natural possibility is to choose a potential $A(x)=A^{(tr)}(x)$ satisfying transversal gauge condition (\ref{eq:TRG}). Note that this condition is fullfilled
  (for $x\neq 0$) for A-B  potential   (\ref{eq:ABp}). Potential $ A^{(tr)}(x)$ decays always as $|x|^{-1}$ at infinity, and hence it can be expected that the SM has a
 structure typical for long-range scattering. However
 for $d=3$ and a given magnetic field $B(x)$, we can also construct a short-range
  magnetic  potential $A(x)$ satisfying equation (\ref{eq:potfield}) and condition
(\ref{eq:H1sr}). Moreover, if $B(x)$ has compact support, then $A(x)$ is also of compact support. This explains why scattering for $d=3$  has a short-range nature.

The difference between dimensions
 $d=3$ and $d=2$ is of topological nature:
 if
$d\geq 3$, then the set ${\Bbb R}^d\setminus\{0\}$ is simply connected whereas it is not true for $d=2$.  Gauge transformations cannot change the magnetic flux if $d=2$. On the contrary, the absence of  a similar invariant for $d=3$ makes the three-dimensional problem essentially more flexible. If $d=2$, then a magnetic potential   cannot even satisfy the condition $A (x)=
o(|x|^{-1}) $ as $|x|\ri\infty$. Indeed, in this case  the second integral in (\ref{eq:flux})  (the circulation of $A(x)$ over the circle
$|x|=R$) tends to $0$ as $R\ri\infty$, and equality (\ref{eq:flux}) implies that 
necessarily $\Phi= 0$.

As was already noted,    condition (\ref{eq:B1}) is
precise, that is the long-range  A-B  effect   occurs even in dimension three if $B(x)$ decays as $|x|^{-2}$ only.  Thus, a long-range behaviour   
 of a magnetic field   in dimension $3$
plays the same role as the topological obstruction 
  in dimension $2$ if the flux $\Phi\not\in
2\pi{\Bbb Z}$.
Of course, our results for $d=3$ remain true for all dimensions $d \geq 3$
(in the general case
 one has to consider $A$ as a 1-form and $B=dA$ as a 2-form).
 Electric potentials are supposed to be zero since they do not add anything new to the phenomena discussed here.

In the next section we discuss some elementary facts  about pseudodifferential operators (PDO) acting on a manifold (the unit sphere).
Then we recall in Section~3 some  basic results  of scattering theory and discuss the 
behaviour of the SM with respect to gauge transformations. The existence of 
the long-range  A-B effect  is   established in Section~4.  On the contrary, in Section~5 we prove its absence for $d=3$ provided condition (\ref{eq:B1}) is satisfied.  

To a large extent, this paper can be considered as a survey article although it contains also some new results. Moreover, compared to \cite{RY1} and \cite{Y2}, we change the point of view supposing that a magnetic field, rather than a magnetic potential, is given.

\section { Pseudodifferential operators on the unit sphere} 

{\bf 2.1.}
The definition of a PDO  $P$   on
the unit sphere $ {\Bbb S}^{d-1}\subset{\Bbb R}^{d}$ reduces to that on a domain $\Sigma\subset {\Bbb R}^{d-1}$  (see, e.g.,  
\cite{Sh}). Roughly speaking,   for a neighbourhood $\Omega$ of an arbitrary point  $\omega_0\in {\Bbb S}^{d-1}$ and a  diffeomorphism $\varkappa: \Omega\to \Sigma$, an operator  $ P : C_{0}^\infty ( \Omega)\to C ^\infty (\Omega)$ reduces by the change of variables $\zeta=\varkappa(\omega)$ to an operator
 $ P _{\varkappa}: C_{0}^\infty ( \Sigma)\to C ^\infty (\Sigma)$. 
 Suppose that, for all $\Omega$ and $\varkappa$, the operators $ P _{\varkappa}$ are PDO on $\Sigma$, that is 
\[
 (P_\varkappa u)(\zeta)=(2\pi)^{-(d-1)/2}\int_{{\Bbb R}^{d-1}}e^{i\langle y,\zeta\rangle}
p_\varkappa(\zeta,y) \hat{u}(y)dy ,
\]
where $\hat{u}$ is the Fourier transform of $u$. Then  $P$  is called 	  PDO on $ {\Bbb S}^{d-1} $.
We require that symbols $p_\varkappa \in C^\infty (\Sigma \times {\Bbb R}^{d-1} )$ of the PDO $P_\varkappa$ belong to
 the (H\"ormander) class   ${\cal S}^{m} $ of symbols 
satisfying, for  all multi-indices $\alpha$ and $\beta$, the estimates
\[
 |\partial_y^\alpha \partial_\zeta^\beta p_\varkappa(\zeta,y)|\leq
C_{\alpha,\beta} (1+|y|)^{m- |\alpha| }.
\]
 Then we say that the PDO $P$ is also from the class ${\cal S}^{m} $.
In terms of the standard PDO notation, $\zeta$ plays the role of the space variable
 and the  variable $y$ is the dual one.

Actually, it suffices to consider only special   diffeomorphisms. For any $\omega_0\in {\Bbb S}^{d-1}$,  
 let $ \Pi_{\omega_0}$
   be the hyperplane orthogonal to $\omega_0$, and let 
$\Omega=\Omega(\omega_0,\gamma)\subset\S^{d-1}$ be determined by the condition
$\langle\omega,\omega_0\rangle > \gamma > 0$.
Let $\zeta=\varkappa(\omega)$ be the orthogonal
projection of $\omega\in\Omega$ on $ \Pi_{\omega_0}$;
 in particular, we assume that $ \varkappa(\omega_0)=0$.
 We denote by $\Sigma$ the orthogonal
projection of $\Omega$ on the hyperplane
 $\Pi_{\omega_0}$ and identify points $\omega\in\Omega$ and
 $\zeta=\varkappa(\omega)$.
 The hyperplane
$\Pi_{\omega_0}$  can be identified with ${\Bbb R}^{d-1}$.
 Let us also consider the unitary mapping
$Z_\varkappa~: L_2({\Omega})\to L_2(\Sigma)$ defined by the equality
 \[
\left(Z_\varkappa u\right)(\zeta)= (1-|\zeta|^2)^{-1/4}u(\omega), \quad
\zeta=\varkappa(\omega).
\]
 Note that compared to the  standard definition it is convenient for us to add the factor
 $(1-|\zeta|^2)^{-1/4}$ in order to make the operator $Z_\varkappa$ unitary.
 We suppose that 
 the operator $P_\varkappa=Z_\varkappa P Z_\varkappa^*$
is a PDO on $\Sigma\subset {\Bbb R}^{d-1}$  with symbol $p_\varkappa(\zeta,y)$  
 from the class ${\cal S}^m =
{\cal S}^m (\Sigma\times {\Bbb R}^{d-1})$.  The symbol $p_\varkappa(\zeta,y)$ is invariant with
respect to diffeomorphisms of $\Sigma$ up to terms from the class ${\cal
S}^{m-1} $. This invariant part,
 considered modulo functions from ${\cal
S}^{m-1} $, is called the principal symbol
 of the PDO $P_\varkappa$ and will
be denoted $p_\varkappa^{(pr)}$.
 The principal symbol of the PDO $P $ is
correctly defined (it means that it does not depend on a choice of $\varkappa$) on the cotangent bundle $T^*\S^{d-1}$ of
 $\S^{d-1}$ by the
equality
 \[
p (\omega,z)= p_\varkappa^{(pr)}(\zeta,y),\quad |\omega|=1,
\quad \langle\omega, z\rangle=0,
 \]
where $\zeta= \varkappa(\omega)$ and $z=\,^t\varkappa^\prime(\omega)y$ is the orthogonal
projection of $y$ on the hyperplane
$\Pi_\omega$. 
 Note also that   kernels $g(\omega,\omega^\prime)$ 
and $ g_\varkappa (\zeta,\zeta^\prime )$ of the
operators $P$ and $P_\varkappa$ regarded as integral operators in $L_2(\Omega)$ and $L_2(\Sigma)$,
respectively,
  are related by the equation
 \[
 g(\omega,\omega^\prime)= g_\varkappa (\zeta,\zeta^\prime )
(1-|\zeta|^2)^{1/4}(1-|\zeta^\prime|^2)^{1/4}, \quad
 \omega,\omega^\prime\in\Omega.
 \]
It is required that $g(\omega,\omega^\prime)$ be a $C^\infty$-function away from the diagonal
$\omega=\omega^\prime$.

\bigskip
 
 {\bf 2.2.}
 We need  information on the essential spectrum
 of a PDO with a   homogeneous symbol of order zero.
 Below  a function $f\in C^\infty$ is called asymptotically
homogeneous of degree $k$ if     $f(z)=|z|^k f(\hat{z})$, $ \hat{z}=z |z|^{-1}$, for   
  $|z|\geq 1/2 $. Of course, $1/2$ is chosen here for definiteness.
Actually, only the behaviour of $f(z)$ for large $|z| $ is essential.
 Let us denote by $T_1^* {\Bbb S}^{d-1} \subset T^* {\Bbb
S}^{d-1}$ the set of pairs
$(\omega,z)$ such that $\omega,z\in {\Bbb S}^{d-1}$ and $\langle\omega,z\rangle=0$.

 \begin{proposition}\label{SpHom} 
Let $P$ be a PDO on ${\Bbb S}^{d-1}$ from the class
 ${\cal S}^0 $ with   principal
symbol $p (\omega,z)$ asymptotically homogeneous of degree $0$ $($in the variable $z)$.
 Then the
essential spectrum $\sigma_{ess}(P)$ of the operator $P$
in the space $L_2({\Bbb S}^{d-1})$ coincides with the
 image   of the
function  $p (\omega,z)$ restricted to the set $T_1^*\S^{d-1}$.
\end{proposition}

As is well known,   kernel $g(\omega,\omega')$ of a PDO $P$
 regarded as an integral
operator       can be very singular on the
diagonal $\omega=\omega'$. Let us find this singularity. 
 
 \begin{proposition}\label{SpHomXY} 
 Under the  assumptions of Proposition~$\ref{SpHom}$, the
kernel $g(\omega,\omega^\prime)$ of a PDO $P$ with principal
symbol $p(\omega, z)$
admits the representation
\begin{equation}\label{eq:PrVaY}
g(\omega,\omega^\prime) = p^{(av)}(\omega)\delta(\omega,\omega^\prime)+
{\rm P.V.}  
q (\omega, \omega^\prime-\omega), 
\end{equation}
up to terms of order $O(|\omega-\omega^\prime|^{-d+1+\nu})$ 
for any   $\nu < 1$ if $d=2$ and for     $\nu=1$ if $d\geq 3$.
Here $\delta(\omega,\omega^\prime)$ is
 the Dirac function on the unit sphere,  
 \[
p^{(av)} (\omega)=|{\Bbb S}^{d-2}|^{-1}\int_{{\Bbb S}^{d-2}_\omega}
p (\omega,\psi)d\psi,\quad {\Bbb S}^{d-2}_\omega={\Bbb S}^{d-1}\cap \Pi_\omega,
 \]
 \[
q(\omega,\tau)= (2\pi i)^{-d+1}  (d-2)!
\int_{{\Bbb S}^{d-2}_\omega}(p (\omega,\psi)- p^{(av)}(\omega
))(\langle\psi, \tau\rangle-i0)^{-d+1}d\psi,
 \]
so that, in particular, 
 for  all $\omega\in {\Bbb S}^{d-1} $
\begin{equation}\label{eq:svc}
 \int_{{\Bbb S}^{d-2}_\omega}
q (\omega,\varphi)d\varphi=0.
\end{equation} 
\end{proposition}

 Note that the function  $q(\omega,\omega^\prime-\omega)$ in (\ref{eq:PrVaY}) is
homogeneous of degree
$-d+1$ in $\omega^\prime-\omega $, so that due to condition (\ref{eq:svc}) the integral operator
with this kernel is correctly defined (as a bounded operator in $L_2({\Bbb S}^{d-1})$) in the sense of  principal value.  
 Thus, up to an integral operator with a weak diagonal singularity, $P $ is   
 the sum $P_0 $ of the operator of multiplication by 
$p^{(av)} (\omega)$ and of the singular integral operator. To be explicit,
 \[
(P_0 f)(\omega)=
  p^{(av)}(\omega) f(\omega)+\lim_{\varepsilon\to 0}\int_{ |\omega' - \omega| > \varepsilon}
q (\omega,  \omega' - \omega) f(\omega')d\omega'.
 \]

 We emphasize that a PDO $P$ of order zero is determined by its principal symbol only up to terms
  from the class  ${\cal S}^{-1} $. These operators are compact so that, by the Weyl theorem, the essential spectra of all such PDO $P$ are the same. Similarly, singular parts of kernels of all such operators $P$ are given by the same formula (\ref{eq:PrVaY}) and all remainders are $O(|\omega-\omega^\prime|^{-d+1+\nu})$.

\begin{remark}\label{d2} 
Let $d=2$. Then $T_{1}^* {\Bbb S}$ consists of points  $ (\omega,\omega^{(+)})$ and 
$ (\omega,\omega^{(-)})$  where $\omega\in {\Bbb S}$ is arbitrary and
   $\omega^{(+)} $ and  $\omega^{(-)}=-\omega^{(+)} $ are obtained from $\omega$ by rotation at the angle $ \pi/2$ and $- \pi/2$   in the positive
$($counterclockwise$)$ direction. Integral over ${\Bbb S}^{0}_\omega$
reduces to a sum over two points $\omega^{(+)} $ and $\omega^{(-)} $ and
\[
\langle\omega^{(\pm)}, \omega^\prime-\omega\rangle=  \pm
|\omega' - \omega| {\rm sgn}\{\omega,\omega'\}+ O(|\omega' - \omega|^2) ,
 \quad \omega' \to \omega,
\]
where $\{\omega,\omega'\}$ is the oriented angle between an
 initial vector $\omega$ and a final vector
$\omega'$. Let us set
\begin{equation}\label{eq:dd2}
  \left.\begin{array}{lcl} 
  p^{(av)}(\omega)&=&2^{-1}(p(\omega,\omega^{(+)}) +p(\omega,\omega^{(-)})),
\\
 p^{(s)}(\omega)&= &(2\pi i)^{-1}(p(\omega,\omega^{(+)}) - p(\omega,\omega^{(-)})). 
  \end{array}\right\}
  \end{equation}
Then formula (\ref{eq:PrVaY}) for
 the singular part of   kernel of the operator $P$ can be written   in the   form
  \[
 g (\omega,\omega')= p^{(av)}(\omega)\delta(\omega,\omega') + p^{(s)}(\omega) {\rm
P.V.}|\omega' - \omega|^{-1}{\rm sgn}\{\omega,\omega'\}.
  \]
\end{remark}
 
Proofs of Propositions~\ref{SpHom} and \ref{SpHomXY} 
can be found in \cite{Y2}.

\section { Scattering matrix } 

{\bf 3.1.}
Let us   discuss briefly some basic facts of scattering theory.
Consider the Schr\"odinger operator
\[
 H=(i\nabla+A(x))^2,\quad x\in{\Bbb R}^d, \quad d\geq 2,
\]
with a real magnetic potential $A(x)=(A_1(x),\ldots,A_d(x))$
satisfying   condition (\ref{eq:H1sr}).
 The dimension $d$ is arbitrary in this section.
We do not assume that the function $A(x)$ is differentiable
so that, strictly speaking, $H$ is correctly defined
as a self-adjoint operator  in the space $ L_2({\Bbb R}^d)$
in terms of the corresponding quadratic form. In general,
equality  (\ref{eq:potfield}) should be understood in the sense of distributions.

Let $H_0=-\Delta$ be the ``free" operator.
Under assumption (\ref{eq:H1sr}) the wave operators
\begin{equation}\label{eq:WW1}
   W_\pm=  W_\pm(H,H_0  ) =  \slim_{t \rightarrow \pm \infty}  e^{iHt}   e^{-iH_0t}   
\end{equation}
 exist, are unitary and possess the intertwining property
\[
 HW_\pm=W_\pm H_0.   
\]
 The scattering operator
${\bf S}=W_+^\ast W_-$ is unitary and commutes with $H_0$.
Let the unitary operator $F :L_2({\Bbb R}^d) \rightarrow
L_2({\Bbb R}_+; L_2({\Bbb S}^{d-1}))$ be defined by the formula
\[
 (F f)(\lambda;\omega)=2^{-1/2}\lambda^{(d-2)/4}\hat{f}(\lambda^{1/2}\omega),
\quad  \lambda> 0, \quad \omega\in {\Bbb S}^{d-1},
\]
where $\hat{f}={\cal F} f$ is the Fourier transform of $f\in L_2({\Bbb R}^d) $.
Clearly, $(F H_{0} f)(\lambda)= \lambda (Ff)(\lambda)$. 
Since ${\bf S}H_0=H_0 {\bf S}$, we have that
\[
(F{\bf S} f)(\lambda)=S(\lambda)(Ff)(\lambda)
\]
where the unitary operator
$S(\lambda): L_2({\Bbb S}^{d-1})\rightarrow L_2({\Bbb S}^{d-1})$
is  known as the scattering matrix (SM).
The scattering amplitude $s(\omega,\omega^\prime;\lambda)$,
$\omega,\omega^\prime\in{\Bbb S}^{d-1}$, is kernel of
$S(\lambda)$ regarded as integral operator.
 The following assertion is well known (see, e.g.,
\cite{Y5, Y2}).

\begin{proposition}\label{reg1A}
 Let condition $(\ref{eq:H1sr})$   hold. Then  
 the operator $T(\lambda)=S(\lambda)-I$ is compact and it belongs to the trace class
if $ \rho> d$. If $ \rho > d+n$,
$n=0,1,2,\ldots$, then
$T(\lambda) $ is integral operator with kernel from the class
$C^n({\Bbb S}^{d-1}\times {\Bbb S}^{d-1})$. 

 Let the condition
\begin{equation}\label{eq:WW3}
 |\partial^\alpha  A(x)|  \leq C_\alpha  (1+|x|)^{- \rho-|\alpha| } 
\end{equation} 
  hold for some $\rho \in (1,d)$ and all multi-indices $\alpha$. Then
 the   operator $T(\lambda) $  has integral    kernel which is a
$C^\infty$-function away from the diagonal
$\omega=\omega^\prime$ and   is bounded by $C(\lambda) |\omega-\omega^\prime|^{-d+ \rho }$
 as $\omega^\prime\rightarrow\omega$.  
\end{proposition}

\bigskip

{\bf 3.2.}
We consider also a class of long-range magnetic potentials satisfying the following

\begin{assumption}\label{ass}
 A magnetic potential $A\in C^\infty$ and
\begin{equation} \label{eq:B3}
 A (x)=A^{(\infty)}(x)+A^{(reg)}(x),
\end{equation}
 where $ A^{(\infty)}\in C^\infty({\Bbb R}^d \setminus\{0\})$
is a homogeneous function of degree $-1$  satisfying the transversal condition 
    \begin{equation}\label{eq:WWX}
  \langle A^{(\infty)}(x),x\rangle=0, \quad x\neq 0,
\end{equation}
    and,   for all $\alpha$,
\begin{equation} \label{eq:H3XXZ} 
 |\partial^\alpha A^{(reg)}(x)|=
O(|x|^{-{\rho}-|\alpha|}),\quad  {\rho}>1.
\end{equation}
\end{assumption}
    
    It turns out that for such long-range potentials the usual wave operators exist. This fact was first observed in
  \cite{LOSTHA}; see also \cite{RY1}, for a different approach.
Nevertheless the structures of the SM are completely different in the short- and long-range cases.  In the long-range case it is natural to regard the SM as a  PDO  on the unit sphere. Its principal symbol can be expressed in terms of  the circulation
\begin{equation}\label{eq:EABAB}
I(x,\xi)=   \int_{-\infty}^\infty
  \langle A ^{(\infty)}(x+ t\xi),\xi\rangle   dt,\quad x\neq 0,\quad \xi\neq 0,\quad \langle x, \xi \rangle =0,  
\end{equation}
 of the homogeneous part 
$A^{(\infty)} $ of $A $ over the straight line $x+ t\xi$ where $t$ runs over ${\Bbb R}$.
 It follows from condition (\ref{eq:WWX})   that
$$\langle A^{(\infty)} (x+ t\xi),\xi\rangle =-t^{-1}\langle A^{(\infty)}  (x+ t\xi),x\rangle=O(|t|^{-2})$$
  as $|t| \rightarrow\infty$, 
 and hence integral  $(\ref{eq:EABAB}) $ converges.
Making in (\ref{eq:EABAB}) the change of variables $t=|x||\xi|^{-1}s$, we arrive at the identity
\begin{equation}\label{eq:Thet}
I(x,\xi)=  I(\hat{x},\hat{\xi}),\quad \hat{x}=x/|x|,\quad \hat{\xi}=\xi/|\xi|,
\end{equation}
 i.e. the function
$I (x,\xi)$ is  
 homogeneous  of degree $0$ in both variables.
Note also that 
\begin{equation}\label{eq:symm}
I(x,-\xi)=-I(x,\xi).
 \end{equation}

\begin{theorem}\label{SM4}
 Let Assumption~\ref{ass}  be satisfied, and let  $S_0$ be the PDO from the class ${\cal S}^0 $
with principal symbol   
\begin{equation}\label{eq:SMS}
 p(\omega, z )=\eta(z)\exp\Bigl(i I(-  z, \omega)\Bigr),\quad
\omega\in{\Bbb S}^{d-1},\; z \in{\Bbb R}^d,\; \langle\omega, z \rangle=0,
\end{equation}
 where  $I(x,\xi)$ is integral $(\ref{eq:EABAB}) $  and  $\eta\in C^\infty$, $\eta(z)=0$ near zero, $\eta(z)=1$ for   $|z|\geq 1/2$.
Then wave operators $(\ref{eq:WW1}) $ exist and  the corresponding SM  admits, for every $p$, the representation
\begin{equation}
 S(\lambda )=S_0+S_p(\lambda)+ \tilde{S}_p(\lambda),
\label{eq:SMpr}\end{equation}
where $S_p(\lambda)$ is a PDO from the class ${\cal S}^{-\nu} $, $\nu=\min\{\rho-1,1\}$,
and   kernel of the operator $ \tilde{S}_p(\lambda)$ is a $C^p$-function of 
$\omega,\omega^\prime\in{\Bbb S}^{d-1}$. 
 \end{theorem}

This  result  almost coincides with Theorem~5.2 of \cite{Y2}. On the other hand, it is a very
 particular case of the general result  of \cite{Y1} where a complete description
of the   amplitude of the PDO $S(\lambda )$
  was obtained for all potentials satisfying condition (\ref{eq:WW3}) for some $\rho > 0$.   Theorem~\ref{SM4} is specially adapted to magnetic potentials $A (x)$ arising naturally from magnetic fields. The operator $S_{0}$ can be considered as the first Born approximation to the SM. Of course the PDO $S_{0}$
 is not determined uniquely by its principal symbol (\ref{eq:SMS}), but the difference of two PDO with the same principal symbol can be included in the operator $S_{p}(\lambda)$.

\bigskip

{\bf 3.3.}
Let us now discuss the behaviour of the SM with respect to gauge transformations
defined by the formula
\begin{equation}\label{eq:GT1}
 \tilde{H}=e^{i\phi}H e^{-i\phi}=(i\nabla+\tilde{A}(x))^2  
\end{equation}
where  
\begin{equation}\label{eq:GT}
 \tilde{A}(x)=A(x)+{\rm grad}\,\phi(x).  
\end{equation}
Of course ${\rm curl}\,\tilde{A}(x)={\rm curl}\,A(x) $.
We are particularly interested in functions $\phi(x) $    which are
asymptotically homogeneous    of degree
zero.  
 
Let us find a relation between    the wave operators
$ W(H,H_0)$ and $ W(\tilde{H},H_0)$.

\begin{proposition}\label{GIW}
Let the wave operators
$W_\pm(H,H_0)$ exist, and
let a differentiable function $\phi  $ be 
be such that 
  $\phi(x)=\phi_0( x)+ \phi_1(x)$ where $\phi_0(  x)=\phi_0(\hat{x})$
 and $ \phi_1(x)=o(1)$ as $ |x|\rightarrow\infty$.
  Then the wave operators $ W_\pm(\tilde{H},H_0)$
also exist and
\begin{equation}\label{eq:GT2}
 W_\pm(\tilde{H},H_0)=e^{i\phi (x)}W_\pm(H,H_0){\cal F}^*
 e^{-i\phi_0(\pm\xi)}{\cal F}.
\end{equation}
\end{proposition}

{\it Proof.} --  Since
$$ 
(\exp(-iH_0t)f)(x)=e^{i|x|^2/(4t)}(2it)^{-d/2}\hat{f}(x/(2t))+o(1),
$$ 
we have that
\begin{eqnarray*}
  (e^{-i\phi}\exp(-iH_0t)f)(x)&=&e^{i|x|^2/(4t)}(2it)^{-d/2} e^{-i\phi_0(\pm
x/(2t))}
\hat{f}(x/(2t))+o(1)
\\ &=&e^{i|x|^2/(4t)}(2it)^{-d/2}\hat{f}^{(\pm)}(x/(2t))+o(1),\quad t\rightarrow\pm\infty,
\end{eqnarray*}
where
$\hat{f}^{(\pm)}(\xi)=e^{-i\phi_0(\pm\xi)}\hat{f}(\xi)$
and the remainder $o(1)$ tends to $0$ in $L_2({\Bbb R}^d)$
 as $ t\rightarrow \pm\infty$. 
 This is equivalent to the relation
\[
 e^{-i\phi}\exp(-iH_0t)f=\exp(-iH_0t)f^{(\pm)}+o(1),
\]
which, in view of 
definition (\ref{eq:GT1}), implies that
\begin{eqnarray*}
 W_\pm(\tilde{H},H_0)f =
\lim_{t\rightarrow\pm\infty} e^{i\tilde{H}t}  e^{-iH_0t}f=
\lim_{t\rightarrow\pm\infty} e^{i\phi}e^{iHt} e^{-i\phi}e^{-iH_0t}f
 \\
 = \lim_{t\rightarrow\pm\infty}e^{i\phi} e^{iHt}  e^{-iH_0t}f^{(\pm)}=
 e^{i\phi}W_\pm(H,H_0)f^{(\pm)}.
\end{eqnarray*}
This proves (\ref{eq:GT2}).  $\quad\Box$

As an immediate consequence of Proposition~\ref{GIW}, we obtain
a relation between      the   scattering operators and matrices. 

\begin{proposition}\label{GI}
Under the assumptions of Proposition~$\ref{GIW}$,
the scattering operators  
are related by the equation
\begin{equation}\label{eq:GT3}
 {\cal F}{\bf S}(\tilde{H},H_0){\cal F}^* = e^{i\phi_0(\xi)}{\cal
F}{\bf S}(H,H_0) {\cal F}^* e^{-i\phi_0(-\xi)}.
\end{equation}
The corresponding SM  $S(\lambda)=S(H,H_0;\lambda)$ and
  $\tilde{S}(\lambda)=S(\tilde{H},H_0;\lambda)$ satisfy
 for all $\lambda > 0$ the relations
\begin{equation}\label{eq:GT1S}
 \tilde{S}(\lambda) =e^{i\phi_0(\omega) }
S(\lambda)e^{ -i\phi_0(-\omega)}
\end{equation}
or, in terms of the scattering amplitudes,
\begin{equation}\label{eq:GT1SA}
 \tilde{s}(\omega,\omega';\lambda)=e^{i\phi_0(\omega)-i\phi_0(-\omega')} s(\omega,\omega';\lambda).
\end{equation}
\end{proposition}

We emphasize that relations (\ref{eq:GT3}) -- (\ref{eq:GT1SA}) for
the scattering operators and
 matrices (but not (\ref{eq:GT2}) for the wave operators)  depend only on the 
asymptotics $\phi_0$ of the  function $\phi$. Probably, formulas (\ref{eq:GT2}),
(\ref{eq:GT3}) appeared first in the paper \cite{RUI} in the case $d=2$ under
some assumptions on $A$. Actually, these formulas do not require any assumptions at all.
 
 As an example, let us consider the case $B(x)=0$.
 
 \begin{example}\label{zero}
Let 
 \[
  A(x)={\rm grad}\,\phi(x),  
\]
where $\phi(x)$ satisfies  the assumptions of Proposition~$\ref{GIW}$.
 Then the wave operators $ W_\pm(H,H_0)$  exist and
\[
 W_\pm(H,H_0)=e^{i\phi (x)} {\cal F}^*
 e^{-i\phi_0(\pm\xi)}{\cal F},
\]
 \[
  {\bf S}(H,H_0)  = {\cal F}^*  e^{i\phi_0(\xi) -i\phi_0(-\xi)}  {\cal F} .
\]
 The corresponding SM does not depend on $\lambda$ and
 \[
 S(H,H_0)   =e^{i\phi_0(\omega) -i\phi_0(-\omega)}.
\]
\end{example}

Relations (\ref{eq:GT1S}) or (\ref{eq:GT1SA}) show  that the SM is not determined by  the magnetic
field $B(x)={\rm curl}\,A(x) $ only although we have an explicit formula which connects the SM in
different gauges. This seems to contradict the following mental experiment. Suppose that a quantum
particle interacts with a magnetic field.
 Note that it is exactly a field but not a potential
 which can be created by our hands. However, to
 calculate the SM theoretically, we have to introduce a magnetic potential
and then solve the Schr\"odinger equation. Thus, the SM
depends on a potential. So it
appears that a particle itself
 chooses a gauge convenient for it. There could be (at least) two possible explanations of this
seeming contradiction. The first is that the scattering amplitude
$s(\omega,\omega^\prime;\lambda)$ cannot be measured
 experimentally although it is widely believed to be possible. From this point of view only the  (differential) scattering cross section
\begin{equation}\label{eq:cross}
\Sigma_{diff}(\omega ,\omega^\prime,\lambda)=(2\pi)^{d-1}
 \lambda^{-(d-1)/2} |  s (\omega ,\omega^\prime;\lambda)|^2, \quad\omega\neq\omega^\prime,
\end{equation} 
($\omega^\prime$ is  an incident direction   of a beam of particles and 
$\omega$ is a  direction of observation) can be practically  found which is
compatible with (\ref{eq:GT1SA}). Another point of view is that 
 experimental devices used for observation of a quantum particle  are not harmless and fix some specific gauge.

On the other hand, for a given field, 
the SM is stable with respect to short-range perturbations
of a potential.

\begin{proposition}\label{SMsr}
Let the wave operators
$W_\pm(H,H_0)$ exist.
Suppose that $ {\rm curl}\, A (x)= {\rm curl}\, \tilde{A} (x)$ and
\begin{equation}\label{eq:SMsr1}
   \tilde{A} (x)-A (x)=O(|x|^{-\rho}), \quad \rho > 1,
\end{equation}
as $|x|\rightarrow \infty$.
Then the wave operators $ W_\pm(\tilde{H},H_0)$
also exist and the scattering operators and matrices
for the pairs $H_0 ,H $ and $H_0 ,\tilde{H} $
coincide.
\end{proposition}

{\it Proof.} --
 According to Propositions~\ref{GIW} and \ref{GI}, it suffices to show
that $A (x)$ and $ \tilde{A} (x)$ are related by equality (\ref{eq:GT})
where the function $\phi(x)$ has a limit (which does not depend on $\hat{x}$) as $|x|\rightarrow \infty$.
Let us define $\phi(x)$     as a curvilinear integral
\[
 \phi(x)= \int_{\Gamma_x} \langle\tilde{A} (y)- A (y), dy\rangle
\]
 taken between  $0$ and a variable point $x$. By the Stokes theorem, 
this integral does not depend on a
choice of $\Gamma_x$ which implies that equality (\ref{eq:GT}) holds. 
 Moreover, choosing $\Gamma_x$ as the piece of straight line connecting $0$ and $x=r\omega$,
$ \omega\in{\Bbb S}^{d-1}$, and using (\ref{eq:SMsr1}), we see that the limit of
$\phi(r\omega)$ as $r \rightarrow \infty$ exists. It remains to show that this limit does not
depend on $ \omega\in{\Bbb S}^{d-1}$. Again by the Stokes theorem,
\begin{equation}\label{eq:CurvB}
 \phi(r\omega_2) - \phi(r\omega_1)= 
\int_{{\Bbb S}_r (\omega_1,\omega_2)} \langle\tilde{A} (y)- A (y), dy\rangle,
\end{equation} 
where ${\Bbb S}_r (\omega_1,\omega_2) $ is the arc of the circle
 centered at the origin and passing
through the points $ r\omega_1$ and $ r\omega_2$.
Condition (\ref{eq:SMsr1}) implies that integral (\ref{eq:CurvB}) 
tends to $0$ as  $r \rightarrow \infty$.$\quad\Box$
 
 \section{Long-range Aharonov-Bohm effect}

 {\bf 4.1.}
 Let us first discuss the case $d=2$. For a given
    magnetic field $B(x)=(0,0,B(x))$, $ x\in {\Bbb R}^2$,  the magnetic
 potential $A^{(tr)}(x)=(A^{(tr)}_1(x),A^{(tr)}_2(x),0)$ satisfying equation
  (\ref{eq:potfield}) (or (\ref{eq:potfield2})) and obeying the transversal gauge condition
 (\ref{eq:TRG})
    can be constructed by the   formulas
 \begin{equation} \label{eq:B2}
 A_1^{(tr)}(x)=-x_2 \int_0^1 B(sx)s ds,
\quad A_2^{(tr)}(x)= x_1 \int_0^1 B(sx)s ds.
\end{equation}
If  condition  (\ref{eq:B1}) is satisfied, then  
 it follows from
  (\ref{eq:B2}) that $A^{(tr)}(x)$ admits representation 
(\ref{eq:B3})
 where $ A^{(\infty)}$
is a homogeneous function of degree $-1$ and $ A^{(reg)} (x)=O(|x|^{-\rho})$ with   ${\rho}= r -1 $ as $|x| \rightarrow\infty$. Indeed, $A^{(\infty)}$ is
 given   by the formula 
 \begin{equation}\label{eq:A'}
 A^{(\infty)}(x)=a(\hat{x})(-x_2,x_1,0) |x|^{-2},\quad \hat{x}=x/|x|,
\end{equation}
where
 \begin{equation} \label{eq:B4bb}
  a(\hat{x})= \int_0^\infty 
B(s\hat{x})s ds 
\end{equation}
  is a   function on the unit circle  and
   \[
   A^{(reg)} (x)=|x|^{-2} (x_{2}, -x_{1})\int_{|x|}^\infty B(s\hat{x})s ds.
   \]
 Moreover,  if $B \in C^\infty({\Bbb R}^d)$ and satisfies the condition
\begin{equation} \label{eq:B1x}
 |\partial^\alpha B(x)|
\leq C_\alpha (1+|x|)^{-r-|\alpha|},\quad r > 2, \quad\forall\alpha,
\end{equation}
 then  $ A^{(tr)}\in C^\infty({\Bbb R}^d)$, $ a\in C^\infty({\Bbb S} )$ and estimates
(\ref{eq:H3XXZ}) hold   for all $\alpha$.
  
    Since ${\rm curl}\,A^{(tr)}(x)=O(|x|^{-r})$  and ${\rm curl}\,A^{(reg)}(x)=O(|x|^{-r})$, it follows from representation (\ref{eq:B3}) that ${\rm curl}\,A^{(\infty)}(x)=O(|x|^{-r})$ where   $r > 2$. On the other hand, ${\rm curl}\, A^{(\infty)}(x)$ is a homogeneous function of degree $-2$ so that necessarily
\begin{equation}\label{eq:reg1}
{\rm curl}\,A^{(\infty)}(x)=0, \quad x\neq 0.
\end{equation}
    The same arguments (or  representation (\ref{eq:A'})) show 
  that the transversal condition (\ref{eq:WWX}) is satisfied.
  Thus, the potential    $A^{(tr)}(x)$ satisfies Assumption~\ref{ass}.

In view of the equalities  (\ref{eq:flux}) and (\ref{eq:A'}),  the total magnetic flux equals
 \begin{equation}\label{eq:flux2}
\Phi  = \int_{\S} a(\psi)d \psi.
\end{equation}
Recall that $\omega^{(\pm)}$ is obtained from $\omega\in{\Bbb S}$ by rotation at the 
angle $\pm \pi/2$ in the positive (counter-clockwise) direction.
 Set
\begin{equation}
 f(\omega)= \int_{\S(\omega^{(-)},\omega^{(+)})} a(\psi)d\psi,
\quad \omega\in{\Bbb S},
\label{eq:B5}\end{equation}
where the integral is taken in the positive  
direction over the half-circle between the points $\omega^{(-)}$ and
$\omega^{(+)}$. Then for any $\omega\in{\Bbb S}$
 \begin{equation}
 f(\omega)+ f(-\omega) =\Phi.
\label{eq:B5fl}\end{equation}
 Comparing formulas (\ref{eq:B4bb}) and  (\ref{eq:B5}),
 we can express the function
$f(\omega)$ in terms of the magnetic field
\begin{equation}
 f(\omega)=  \int_{\langle x,\omega\rangle\geq 0}B(x)dx.
\label{eq:B5B}\end{equation}
  In its turn,  integral (\ref{eq:EABAB}) can be expressed   in terms of the
   function $f(\omega)$.
  
  \begin{lemma}\label{A1}
For all $\omega\in {\Bbb S}^{d-1}$, we have that  
\begin{equation}
 I (\omega,\omega^{(\pm)})= \pm f (  \omega).
\label{eq:Thg}
\end{equation}
\end{lemma}
  
{\it Proof.} --
 By virtue of (\ref{eq:symm}), it suffices to consider the case of the upper
sign. Since
\[
<(-\omega_2-t\omega_2^{(+)}, \omega_1+t\omega_1^{(+)}), (\omega_1^{(+)},\omega_2^{(+)})>
= \omega_1 \omega_2^{(+)} -\omega_2 \omega_1^{(+)}=1,
\]
we have that for potentials (\ref{eq:A'})
$$
I (\omega,\omega^{(+)})=\int_{-\infty}^\infty
a \Bigl(\frac{\omega+t\omega^{(+)}}{\sqrt{t^2+1}}\Bigr)\frac{dt}{t^2+1}.
$$
 Making the change of variables $t=\tan \psi$, we
get formula (\ref{eq:Thg}).$\quad\Box$

We denote by $ S^{(tr)}(\lambda)$ the SM corresponding to the potential $ A^{(tr)}$.
Given Theorem~\ref{SM4} the following two assertions are immediate consequences of Propositions~\ref{SpHom}
and \ref{SpHomXY} (see also Remark~\ref{d2}).

\begin{theorem}\label{ABStil}
Let $d=2$ and let condition $(\ref{eq:B1x})$ be satisfied.
Define the function $f(\omega)$   by formula $(\ref{eq:B5B})$ 
and set $\gamma_+= \max f(\omega)$,  $\gamma_-= \min f(\omega)$.
 Then for all $\lambda > 0$ relation 
\begin{equation}\label{eq:AhBo}
\sigma_{ess}(S^{(tr)}(\lambda))=[\exp(i\gamma_-), \exp(i\gamma_+)]\cup
[\exp(-i\gamma_+), \exp(-i\gamma_-)]
\end{equation}
 holds     if $\gamma_+-\gamma_- < 2\pi$,
and   $\sigma_{ess}(S^{(tr)}(\lambda))$ covers the whole unit circle ${\Bbb T}\subset {\Bbb C}$
 if $\gamma_+-\gamma_-\geq 2\pi$.
\end{theorem}

{\it Proof.} --
Let us apply Proposition~$\ref{SpHom}$ to the PDO $S_{0}$ with principal symbol (\ref{eq:SMS}).
It follows from formula (\ref{eq:Thg}) that in this case
\begin{equation}\label{eq:AhBoXY}
p(\omega,\omega^{(\pm)})=e^{\pm i f(\omega^{(\mp)})}.
\end{equation}
Therefore the images of the functions $p(\omega,\omega^{(+)})$ and $p(\omega,\omega^{(-)})$ coincide  with the first and the second arcs 
in (\ref{eq:AhBo}), respectively, if $\gamma_+-\gamma_- < 2\pi$.
If $\gamma_+-\gamma_-\geq 2\pi$, then each of these images covers the whole unit circle. So it remains to take into account that according to representation
(\ref{eq:SMpr}) the essential spectra of the operators $S^{(tr)}(\lambda)$ and $S_{0} $ are the same.
$\quad\Box$ 

\begin{theorem}\label{ABStilde}
Let the assumptions of Theorem~$\ref{ABStil}$ hold.
  Set 
\begin{equation}\label{eq:SAB}
s_0 (\omega,\omega') = e^{i (f(\omega^{(-)})-f(\omega^{(+)}))/2}
\Bigl( \cos(\Phi/2)\delta(\omega,\omega')
+\pi^{-1}\sin(\Phi/2)
{\rm P.V.}\frac{{\rm sgn}\{\omega,\omega'\}}{|\omega - \omega'|}\Bigr).
\end{equation}
  Then for arbitrary $\lambda >0$ 
\[
 |s^{(tr)} (\omega,\omega';\lambda)-s_{0}
(\omega,\omega')|\leq C(\lambda) |\omega-\omega^\prime|^{- 3+r_{0}}.
\]
Here $r_{0}= r   $ if $r  < 3$ and  $r_{0} $
is an arbitrary number smaller than $3$ if $r\geq 3$.
\end{theorem}

{\it Proof.} --
Now we apply Proposition~$\ref{SpHomXY}$ to the PDO $S_{0}$ with principal symbol (\ref{eq:SMS}).
Comparing formulas \e{eq:dd2} and \e{eq:AhBoXY}, we see that in the case considered
\[
p^{(av)}(\omega)=2^{-1}(e^{ i f(\omega^{(-)})}+ e^{-i f(\omega^{(+)})}),
\q
p^{(s)}(\omega)= (2\pi i)^{-1}(e^{ i f(\omega^{(-)})}- e^{-i f(\omega^{(+)})}).
\]
Using identity \e{eq:B5fl}, we find that formula  \e{eq:dd2} yields expression \e{eq:SAB} for the singular part of kernel of operator $S_{0}$. The ``regular" part of its kernel       is $O( |\omega-\omega^\prime|^{-  \varepsilon})$ as $|\omega-\omega^\prime| \to 0$ for any $\varepsilon>0$. Kernel of the operator $S_{p}(\lambda)$ in \e{eq:SMpr} satisfies the same  estimate if $r\geq 3$ and it is $O( |\omega-\omega^\prime|^{ -3+r})$ if $r<3$.
$\quad\Box$ 

\begin{corollary}\label{ABStildeCR}
The diagonal singularity of the scattering cross section $(\ref{eq:cross})$ is given by the formula
\begin{equation}\label{eq:Sigm}
\Sigma_{diff}(\omega ,\omega^\prime;\lambda)=2\pi^{-1}
 \lambda^{ -1/2} \sin^2 (\Phi/2)\, |\omega-\omega^\prime|^{-2} +O(|\omega-\omega^\prime|^{- 4+ r_{0}}).
\end{equation}
\end{corollary}

Thus, the singular part $S_0$ of the SM $S^{(tr)}(\lambda)$ is
the integral operator in $L_2({\Bbb S})$ with kernel (\ref{eq:SAB}).
Up to the phase factor, it is determined by the magnetic flux $\Phi$ only (and does not depend on $\lambda$).
 We see that  in the dimension two
 even for magnetic fields of
compact support with $\Phi\not\in 2\pi {\Bbb Z}$,
 the SM contains the singular integral operator
and the forward singularity 
(\ref{eq:Sigm}) of the
scattering cross section  
 is stronger than for short-range
magnetic potentials where it is $O(|\omega-\omega^\prime|^{-4+2\rho})$.
On the contrary, if
   $\Phi\in 2\pi{\Bbb Z}$, then  according to (\ref{eq:B5fl}), (\ref{eq:SAB})
the operator $S_{0}$ acts as multiplication
by the function $e^{i f(\omega^{(-)})}$.
  As we shall see  in the next section,
this situation is typical for dimensions $d\geq 3$.
Note also that if $B(x)$ is an even function, that is $B(x)=B(-x)$, then, again by (\ref{eq:B5fl}),
$f(\omega)=\Phi/2$ for all $\omega\in{\Bbb S}$,
 and hence the first factor in the right-hand side
of (\ref{eq:SAB}) equals $1$. In this case $\sigma_{ess}(S(\lambda))$ consists
of the two points $e^{i\Phi/2}$ and $e^{-i\Phi/2}$. Of course in the case $a(\omega)=-\alpha$ formula (\ref{eq:SAB}) coincides, up to smooth terms, with formula (\ref{eq:kern}) if the natural parametrization of the unit circle ${\Bbb S}$ by points $\theta\in  [0,2\pi)$ is used.

As a concrete example, let us consider the field
 \[
 B(x)=B_{0}(r)+ B_{1}(r) \langle q,\hat{x} \rangle, \quad r =|x|,
 \]
 where $B_{0}$ and $ B_{1}$ are $C^\infty$-functions with compact supports
 and $q\in{\Bbb R}^2$ is some given vector.
It follows from formula  (\ref{eq:B4bb}) that
in this case
$$ 
a(\hat{x})=-\alpha +  \langle p,\hat{x}\rangle,\quad \alpha\in{\Bbb R},\quad p\in{\Bbb R}^2,
$$
where
\[
\alpha=-\int_{0}^\infty B_{0}(r) rdr, \quad
p= q \int_{0}^\infty B_{1}(r) rdr .
\]
 Clearly, $\Phi=-2\pi\alpha$. Let us calculate   function (\ref{eq:B5}). For an arbitrary
$\omega\in {\Bbb S}$, let $\varphi$ be the angle between $\omega$ and $p$, and let $\theta$ be the angle between $\hat{x}$ and $p$. Then $a(\hat{x})=-\alpha+|p| \cos\theta$ and
$$ 
f(\omega)=-\pi\alpha+ |p| \int _{\varphi-\pi/2}^{\varphi+\pi/2}\cos\theta d\theta=-\pi\alpha +2| p|\cos\varphi=-\pi\alpha+ 2\langle p,\omega\rangle.
$$ 
Therefore the conclusion of Theorem~\ref{ABStil} is true with  
$\gamma_+=  - \pi\alpha+2|p|  $ and $\gamma_-= - \pi\alpha -2|p| $.
 In particular,  if  $2|p|\geq \pi$, then $\sigma_{ess}(S^{(tr)}(\lambda))={\Bbb T}$.  On the
contrary,  $\sigma_{ess}(S^{(tr)}(\lambda))$ consists of the two points
$\exp(\pi i \alpha)$ and $\exp(-\pi i \alpha)$ if $p=0$.
The phase factor in (\ref{eq:SAB}) equals $ \exp(2 i \langle p,\omega^{(-)}\rangle)$.

Actually, the results above do not require that the potential satisfy transversal condition
(\ref{eq:TRG}). The next result follows again from
  Theorem~\ref{SM4} combined with Propositions~\ref{SpHom}
and \ref{SpHomXY}.

\begin{theorem}\label{ABSt}
Suppose that $A\in C^\infty({\Bbb R}^2)$ admits representation $(\ref{eq:B3})$ where $A^{(\infty)}$ is function $(\ref{eq:A'})$ with $a \in C^\infty({\Bbb S})$ and
$A^{(reg)}$ satisfies estimates $(\ref{eq:H3XXZ})$. Let $f$ be function 
$(\ref{eq:B5})$. Then all conclusions of Theorems~$\ref{ABStil}$ and $\ref{ABStilde}$
remain true for the SM corresponding to the potential $A$.
\end{theorem}
  
Let us discuss this result from the point of view of gauge transformations.
Let two potentials $A(x)$ and $\tilde{A}(x)$ satisfy the assumptions of Theorem~\ref{ABSt},
and let $\Phi$ and $\tilde{\Phi}$ be  the   corresponding magnetic fluxes.
If they are related by equality (\ref{eq:GT}) where $\phi(x)$ satisfies the conditions of Proposition~\ref{GIW}, then 
$\tilde{\Phi}=\Phi$,
\begin{equation}\label{eq:aaw} 
\tilde{a}(\omega)=a(\omega)+\phi_0^\prime(\omega) 
\end{equation}
and hence according to (\ref{eq:B5})
\[
\tilde{f}(\omega^{(-)})-\tilde{f}(\omega^{(+)})=f(\omega^{(-)})-f(\omega^{(+)})
+2( \phi_0 (\omega)-\phi_0 (-\omega)).
\]
It follows from (\ref{eq:SAB}) that singular parts of the corresponding SM are connected by the equality
\begin{equation}\label{eq:SABn}
 \tilde{s}_0 (\omega,\omega') = e^{i\phi_0 (\omega)-i\phi_0 (-\omega)}s_0 (\omega,\omega'),
\end{equation}
which agrees with exact formula (\ref{eq:GT1SA}) for scattering amplitudes.

 Conversely, if $\tilde{\Phi}=\Phi $, then  
  the function
\begin{equation}\label{eq:aaw1} 
\phi_0 (\omega)=\int_{\S(\omega_0 ,\omega )} (\tilde{a}(\psi)-a(\psi))d \psi
\end{equation}
(the point $\omega_0\in\S$ is arbitrary but fixed) 
 is correctly defined on the unit circle and
equality (\ref{eq:aaw}) is   satisfied. Set $\phi(x)=\eta(|x|) \phi_0(\hat{x})$
where $\eta\in C^\infty$, $\eta(r)=0$ in a neighbourhood of zero and
$\eta(r)=1$ for large  $r $. It follows from (\ref{eq:aaw}) that  equality (\ref{eq:GT}) holds, up
to a term $A_{sr}(x)$ satisfying estimates (\ref{eq:H3XXZ}), that is, $\tilde{A}=A+{\rm
grad}\,\phi+A_{sr}$.
 Therefore the SM $S(\lambda)$ and
$\tilde{S}(\lambda)$ for the Schr\"odinger operators with magnetic potentials $A$ and
 $\tilde{A}-A_{sr}$ are
related by equality (\ref{eq:GT1S}) and hence their
singular parts are related by equality (\ref{eq:SABn}).
This implies that if  (\ref{eq:SAB}) is verified for
$\tilde{A}-A_{sr}$, then it is also true for $A$.
Thus, for a given $\Phi $,
it suffices to prove Theorem~\ref{ABStilde} only for one function $a$
satisfying (\ref{eq:flux2}) (but for all short range terms $A^{(reg)}$).
We can choose $ a (\omega)=(2\pi)^{-1}\Phi$, which reduces the proof of 
Theorem~\ref{ABStilde} to the case of a constant function $a$.
In particular, if $\Phi=0$, then the problem reduces to the short range case.
 The same is true with respect to
 Theorem~\ref{ABStil}   for even functions $ a (\omega)$
only. Then   function (\ref{eq:aaw1}) is
also even  so that, by virtue of (\ref{eq:GT1S}), the SM $\tilde{S}(\lambda)$ and
$S(\lambda)$ are unitarily equivalent.

 \bigskip

 {\bf 4.2.}
    Here we consider arbitrary magnetic potentials
 $ A(x) $, $x\in{\Bbb R}^3$,  with Coulomb decay at infinity
satisfying, at least asymptotically, the transversal gauge condition.
For such   potentials, the magnetic field $B(x)={\rm curl}\,A(x)$
 decays, in general, as $|x|^{-2}$ at infinity, 
so that   assumption (\ref{eq:B1}) is
not satisfied. We shall show that in this case   the SM
contains a singular integral operator and hence the long-range A-B  effect occurs.

The next two results extend Theorems~\ref{ABStil} and \ref{ABStilde} to the case $d=3$.
    They  follow again from Theorem~\ref{SM4} and  Propositions~\ref{SpHom} and \ref{SpHomXY}   applied to  the PDO with symbol (\ref{eq:SMS}).
    
    \begin{theorem}\label{SMsingE}
Suppose that $d=3$. Let Assumption~$\ref{ass}$ be satisfied, and let $I(x,\xi)$ be integral \e{eq:EABAB}. Then $\sigma_{ess}(S(\lambda))$ coincides with the image of the function $I(\psi,\omega)$ for all $\psi,\omega\in{\Bbb S}^2$ such that $\langle\psi, \omega\rangle=0$. 
\end{theorem}

\begin{theorem}\label{SMsing}
 Under the assumptions of Theorem~$\ref{SMsingE}$ define
 the functions  
 \begin{equation}\label{eq:qqq1X}
p^{(av)} (\omega)=(2\pi)^{-1}\int_{{\Bbb S}_\omega}
\exp(iI ( \psi, \omega))d\psi,\quad {\Bbb S}_\omega={\Bbb S}^{2}\cap \Pi_\omega,
\end{equation}
 \begin{equation}\label{eq:qqq1A}
q(\omega,\tau)= -(2\pi )^{-2}   
\int_{{\Bbb S}_\omega}(\exp(iI (\psi, \omega)) - p^{(av)}(\omega
))(\langle\psi, \tau\rangle-i0)^{-2}d\psi,
\end{equation}
and set
 \begin{equation}\label{eq:ssing}
s_0(\omega ,\omega^\prime)
=p^{(av)}(\omega)\delta(\omega,\omega^\prime)+
{\rm P.V.}   q(\omega,  \omega^\prime-\omega ).
\end{equation}
Then for an arbitrary $\lambda > 0$ the scattering amplitude satisfies
 the  estimate
\[
 |s (\omega,\omega';\lambda)-s_{0}
(\omega,\omega')|\leq C(\lambda) |\omega-\omega^\prime|^{- 3+\rho_{0}},
\]
where $\rho_{0}=  \rho $ if $ \rho \in (1,2)$ and $\rho_{0}=2$ if $ \rho \geq 2$.
\end{theorem}

\begin{corollary} \label{B6as}
 If $\omega\neq\omega^\prime$ but $\omega-\omega^\prime\rightarrow 0 $, then
\[
s(\omega ,\omega^\prime ;\lambda)= q(\omega,
 \omega^\prime-\omega )  + O(|\omega -\omega^\prime|^{-3+\rho_{0}}).
\]
 \end{corollary}

\begin{corollary} \label{section}
 If $\omega\rightarrow\omega^{\prime}$, then
\[
\Sigma_{diff}(\omega ,\omega^{\prime};\lambda)=(2\pi)^2
 \lambda^{-1}  |q (\omega ,  \omega^{\prime}-\omega  ) |^2
 + O(|\omega -\omega^{\prime}|^{-5+ \rho_{0}}).
\]
 \end{corollary} 

Note that the order of singularity $|\omega  - \omega^{\prime} |^{-4}$
here is the same as for electric Coulomb potentials. We emphasize also that singular part      
\e{eq:ssing} of the SM does not depend on $\lambda$.
 
  It follows from
equality (\ref{eq:qqq1A}) that if 
the  function  $I(x,\omega)$ does not depend on $x$,  then
$q (\omega,\tau)=0$ so that the singular integral operator
  disappears in (\ref{eq:ssing}). We shall see in the next section that   this situation  really occurs if the magnetic field satisfies condition (\ref{eq:B1}).

   \medskip

  Nevertheless $q (\omega,\tau)$ is non-trivial in the general case.
Let us consider (see \cite{Y2}, for details) two concrete examples  
of   potentials $A^{(\infty)}(x)$ homogeneous of degree $-1$ and satisfying 
transversal condition (\ref{eq:WWX}).
The corresponding
    fields ${\rm curl}\, A(x)$   decay  only as $|x|^{-2}$ at infinity.

 We define the first of these  potentials by the equation
\[
 A^{(\infty)}(x)=|x|^{-3}(\alpha_1 x_2 x_3,\alpha_2  x_3 x_1 , \alpha_3 x_1 x_2),
  \quad x=(x_1,x_2,x_3)\in{\Bbb R}^3,
\]
 where $\alpha_j$ are constants and
$\alpha_1+\alpha_2 + \alpha_3 =0$.
 An easy calculation shows    that 
  function (\ref{eq:EABAB}) equals
  \[
I(x,\omega)= 2   |x|^{-2} (\alpha_1\omega_1  x_2 x_3+\alpha_2 \omega_2 x_3 x_1+
\alpha_3\omega_3 x_1 x_2).
  \]
  Since this function depends on $\lambda$, it cannot be expected that $q (\omega,\tau)=0$. 
  
  Actually, the functions
  $p^{(av)} (\omega)$ and $q(\omega,\tau)$ can be calculated explicitly.
   For an arbitrary $\omega =(\omega_1,
\omega_2,\omega_3)$, the coordinates of
an arbitrary point $x=(x_{1},x_{2},x_{3})\in{\Bbb S}_\omega={\Bbb S}^2 \cap \Pi_\omega$
 can be written, for some $\theta\in[0,2\pi)$, as
\begin{equation}\label{eq:xxx}
   \left.\begin{array}{lcl} 
x_1 = -(\omega_1^2+\omega_2^2)^{-1/2}(\omega_2\cos\theta+\omega_1 \omega_3\sin\theta ),
   \\
x_2 = (\omega_1^2+\omega_2^2)^{-1/2}(\omega_1\cos\theta -\omega_2 \omega_3\sin\theta ),
\quad 
 x_3 = (\omega_1^2+\omega_2^2)^{1/2} \sin\theta.
   \end{array}\right\}
\end{equation}
 Set
\[ 
{\cal A} (\omega)=(\omega_1^2+\omega_2^2)^{-1}
\Bigl(4\alpha_3^2 \omega_1^2\omega_2^2\omega_3^2 
 +(\alpha_1(\omega_1^2-\omega_2^2\omega_3^2)
-\alpha_2(\omega_2^2-\omega_1^2\omega_3^2))^2
\Bigr)^{1/2}  
\] 
and define the angle $\theta_0(\omega)$  by the equations
 \[
  \left.\begin{array}{lcl} 
 \sin \theta_0(\omega)=-2\alpha_3  \omega_1 \omega_2 \omega_3 (\omega_1^2+\omega_2^2)^{-1}
{\cal A} (\omega)^{-1},
\\  
\cos\theta_0(\omega)= 
(\alpha_1(\omega_1^2-\omega_2^2\omega_3^2) -\alpha_2(\omega_2^2-\omega_1^2\omega_3^2))
(\omega_1^2+\omega_2^2)^{-1}{\cal A} (\omega)^{-1}.
   \end{array}\right\}
\]
Then
\[
 p^{(av)}(\omega)=(2\pi)^{-1}\int_0^{2\pi}\cos ({\cal
A}(\omega)\sin\theta)d\theta
\]
and  
\[ 
q(\omega ,\tau) = - (2\pi )^{-2}\int_0^{2\pi}\Bigl(  e^{ i{\cal A}(\omega)
\sin(2\theta+\theta_0(\omega))}-
p^{(av)}(\omega)\Bigr)
 (\langle x(\theta),\tau\rangle -i0)^{-2}d\theta 
\]
where $x(\theta)$ is defined by formulas \e{eq:xxx}.
  
   As another example,
 we choose a modification of the A-B  potential
\[
 A^{(\infty)}(x)=-\alpha |x|^{-2}(- x_2, x_1,0).
\]
In this case
  \[
I(x,\omega)=  \pi\alpha   |x|^{-1} ( \omega_1 x_2-\omega_2 x_1)
    =\pi\alpha(1-\omega_3^2)^{1/2} \cos \theta
    \] 
if $x,\omega$ and $\theta$ are related by equalities (\ref{eq:xxx}). Plugging this expression into (\ref{eq:qqq1X}) and (\ref{eq:qqq1A}), we obtain explicit representations for the functions $p^{(av)} (\omega)=p^{(av)} (\omega_{3})$ and $q(\omega,\tau)
    =q(\omega_{3},\tau)$:
    \[
 p^{(av)}(\omega_{3})=(2\pi)^{-1}\int_0^{2\pi}\cos  (\pi\alpha(1-\omega_3^2)^{1/2} \cos \theta ) d\theta,
\]  
\[
q(\omega_{3} ,\tau) = - (2\pi )^{-2}\int_0^{2\pi}
\Bigl(  e^{  i \pi\alpha(1-\omega_3^2)^{1/2} \cos \theta }
- p^{(av)}(\omega)\Bigr)
  (\langle x(\theta),\tau \rangle -i0)^{-2}d\theta.
\]

In  both these examples the   SM contain   singular integral operators
and hence the long-range A-B effect occurs.

\section{There is no long-range Aharonov-Bohm  effect \\ in dimension three}

{\bf 5.1.}
In this section we suppose that the dimension $d=3$.
Our results remain true for all $d\geq 3$ but not for $d=2$.
Let  $B(x)=(B_1(x),B_2(x),B_3(x))$ be  a magnetic field  
such that  ${\rm div}\,B(x)=0$. Recall that a magnetic potential
$A^{(tr)}(x)=(A_1^{(tr)}(x),A_2^{(tr)}(x),A_3^{(tr)}(x))$ satisfying  equation (\ref{eq:potfield}) and the transversal gauge condition (\ref{eq:TRG})
 is  constructed by the formula
\begin{equation} \label{eq:BB2}
 A_1^{(tr)}(x)=\int_0^1\Bigl( B_2(sx)x_3-  B_3(sx) x_2\Bigr)s ds.
\end{equation}
 Expressions for components $ A_2^{(tr)}(x)$ and $ A_3^{(tr)}(x)$ are  obtained by  
 cyclic permutations of indices in (\ref{eq:BB2}).  If
estimate (\ref{eq:B1}) is satisfied, 
  then $A^{(tr)}(x)$ admits the representation (\ref{eq:B3})
 where 
\begin{equation} \label{eq:Binf}
 A_1^{(\infty)}(x)=|x|^{-2}\int_0^\infty\Bigl( B_2(s\hat{x})x_3-  B_3(s\hat{x}) x_2\Bigr)s ds, 
\end{equation}
\begin{equation} \label{eq:Breg}
 A_1^{(reg)}(x)=-|x|^{-2}\int_{|x|}^\infty\Bigl( B_2(s\hat{x})x_3-  B_3(s\hat{x}) x_2\Bigr)s ds.
\end{equation}
Thus, 
$A^{(\infty)}(x)$ is a homogeneous function of degree $-1$
and $      A^{(reg)}(x)=
O(|x| ^{- \rho}) $
  with  $ \rho=r-1> 1$ as $|x|\ri\infty$.
  Quite similarly to the two-dimensional case (see subsection~4.1), we have that 
$ A^{(\infty)}(x)   $  satisfies  equations (\ref{eq:WWX})  and  
(\ref{eq:reg1}).

Given a magnetic field $B(x)$ obeying condition  (\ref{eq:B1}), we construct now
a magnetic potential $A(x)$ satisfying equation (\ref{eq:potfield}) and  estimate  (\ref{eq:H1sr}). 
We proceed from the magnetic potential $A^{(tr)}$ in the transversal gauge.  
   Let $ A ^{(\infty)}$ be function (\ref{eq:Binf}). 
We define     the function
$U(x)$   for   $ x\neq 0$  as a curvilinear integral
\begin{equation}\label{eq:Curv}
 U(x)= \int_{\Gamma_{x_0,x}} \langle A^{(\infty)}(y), dy\rangle
\end{equation} 
 taken between some fixed point $x_0 \neq 0$ and a variable point $x$. It is required that
  $0 \not\in\Gamma_{x_0,x}$,
 so that, in view of (\ref{eq:reg1}) and the Stokes theorem, $U(x)$ does not depend on a choice of a contour
$\Gamma_{x_0,x}$. Here it  is used that the set ${\Bbb R}^3\setminus\{0\}$  (and ${\Bbb
R}^d\setminus\{0\}$ for all $d\geq 3$)   is simply connected. Clearly, 
\begin{equation}\label{eq:GG}
  A^{(\infty)} (x)={\rm grad}\, U(x).
\end{equation}
 Moreover, the function $U(x)$ is  homogeneous  of degree $0$. Indeed, if $x_2=\gamma x_1$,
$\gamma>1$, then we can choose $\Gamma_{x_0,x_2}= \Gamma_{x_0,x_1}\cup
(x_1,x_2)$ where  
 $ (x_1,x_2)$ is the piece of straight line connecting $x_1$ and $x_2$.
If $y \in(x_1,x_2)$, then according to (\ref{eq:WWX})
$\langle A^{(\infty)}(y), d y\rangle=0$.
Hence $U(x_1)=U (x_2)$. We use definition (\ref{eq:Curv})
away from some neighbourhood of the point $x=0$  and
extend  $U(x)$ as a  differentiable
function to all $ {\Bbb R}^3$. For example, we can choose some numbers $R_{2}> R_{1} > 0$ and a function 
$\eta\in C^\infty ({\Bbb R}^3)$ such that $\eta(x)=1$ for $|x|\geq R_{2}$,
$\eta(x)=0$ for $|x|\leq R_{1}$ and then replace $U(x)$ by
$\eta(x) U(x)$. 
Let us now set
\begin{eqnarray}\label{eq:GGx}
A(x) &=& A^{(tr)}(x)-{\rm grad}\, ( \eta(x) U(x))
\nonumber\\
&=& A^{(reg)}(x) + (1-\eta(x)) A^{(\infty)}(x) -  U(x) {\rm grad}\,\eta(x), 
\end{eqnarray}
so that $A(x)=A^{(reg)}(x)$ for $|x|\geq R_{2}$ and $A(x)=A^{(tr)}(x)$
for $|x|\leq R_{1}$.
 Thus, we have the following result.

\begin{proposition}\label{pot} 
Suppose   that  ${\rm div}\,B(x)=0$ and
that condition $(\ref{eq:B1})$
holds. Let the magnetic potential $A(x)$ be defined by formula 
$(\ref{eq:GGx})$ where $A^{(\infty)}(x)$, $A^{(reg)}(x)$ and $U(x)$
are functions $(\ref{eq:Binf})$, $(\ref{eq:Breg})$ and
$(\ref{eq:Curv})$, respectively.
Then   $A(x)$   satisfies equation $(\ref{eq:potfield})$
and estimate  $(\ref{eq:H1sr})$. Moreover, $A(x)$ has compact support
if $B(x)$ has compact support.
\end{proposition}
 
 In the  case of magnetic fields $B(x)$ with compact supports our construction is   close to that of \cite{HeHe}.
By the proof of Proposition~\ref{pot} we could have proceeded from the magnetic potential
$A^{(c)}(x)$ satisfying the Coulomb gauge condition ${\rm div}\, A^{(c)}(x)=0$.
This is however less convenient.  
 
 Suppose that a magnetic field is supported by some ball ${\Bbb B}$,  ${\Bbb B}^{\prime}$ is a slightly larger ball and a direct interaction of quantum particles with this field is excluded  by the Dirichlet boundary condition on $\partial{\Bbb B}^{\prime}$.
Proposition~\ref{pot} shows that we can choose a magnetic potential supported by  ${\Bbb B}^{\prime}$ so that scattering in this case is trivial. On the other hand, if a magnetic field is supported by some torus ${\bf T}$ and  the Dirichlet boundary condition is put on the boundary of a slightly larger torus ${\bf T}^{\prime}$, then the Stokes theorem does not allow us to find a potential supported by ${\bf T}^{\prime}$ (provided the magnetic flux through a section of ${\bf T}$ is not zero). Therefore scattering in this case is non-trivial although it is of short-range nature.

Let us denote by ${\cal A}(B)$ the class of
 magnetic potentials satisfying equation  (\ref{eq:potfield}) 
and estimate  (\ref{eq:H1sr}) for some $\rho > 1$.
 This class is non-empty according to Proposition~\ref{pot}.
If $A\in{\cal A}(B)$, then, for an arbitrary function
$\phi (x)$ such that ${\rm grad}\, \phi (x)=
O(|x|^{- \rho}) $, potential (\ref{eq:GT}) also belongs to this class.
According to Proposition~\ref{SMsr} the SM
for the pair $H_0=-\Delta$, $H=(i\nabla+A(x))^2$ does not depend on the choice of
$A\in{\cal A}(B)$ and, thus, is determined by 
the magnetic field $B(x)$ only. We say that this SM $S(\lambda)=S(\lambda;B)$
is the SM for the  field $B(x)$.

Comparing Propositions~\ref{reg1A} and \ref{pot}, we arrive at the following result.

\begin{theorem}\label{main} 
Let a magnetic field $B(x)$ be such that
${\rm div}\,B(x)=0$ and
  condition $(\ref{eq:B1})$ holds, and let
a magnetic potential $A\in{\cal A}(B)$. Then the wave operators 
for the pair $H_0=-\Delta$, $H=(i\nabla+A(x))^2$
exist, are unitary and  the   SM
$ S(\lambda)  $ for the magnetic field $B(x)$ is a unitary
operator for all $  \lambda > 0  $.
The operator $T(\lambda)=S(\lambda)-I$ is compact and it belongs to the trace class
if $ r > 4$. If $r> 4+n$,
$n=0,1,2,\ldots$, then
$T(\lambda) $ is integral operator with kernel from the class
$C^n({\Bbb S}^2\times {\Bbb S}^2)$. If  condition $(\ref{eq:B1x})$ holds
    for some $r \in (2,4)$ and all multi-indices $\alpha$, then
 the   operator $T(\lambda) $  has integral    kernel which is a
$C^\infty$-function away from the diagonal
$\omega=\omega^\prime$ and   is bounded by $C(\lambda) |\omega-\omega^\prime|^{-4+ r }$
 as $\omega^\prime\rightarrow\omega$.  
\end{theorem}

\begin{corollary} \label{B6}
 If estimate $(\ref{eq:B1})$ is satisfied for $r> 4$, then 
$\Sigma_{diff}(\omega ,\omega^{\prime}; \lambda)$ is a bounded function of
$\omega,\omega^{\prime}\in{\Bbb S}^2$. If condition $(\ref{eq:B1x})$ is satisfied
 for some $r \in (2,4)$ and all multi-indices $\alpha$, then
\[
\Sigma_{diff}(\omega ,\omega^{\prime} ; \lambda)= 
  O(|\omega -\omega^{\prime}|^{-8+2r}) \quad {\rm as} \quad \omega \rightarrow\omega^{\prime}.
\]
 \end{corollary}

Using the first formula (\ref{eq:GGx}) and applying
Proposition~\ref{GI} to the function $\phi(x)=U(x)$, we can also describe the structure of the SM in the trasversal gauge.

\begin{proposition} \label{B6tr}
Suppose that a magnetic field $B(x)$
satisfies the assumptions of Theorem~$\ref{main}$.
Let $ S^{(tr)}(\lambda) $ be the SM for the  pair $H_0=-\Delta$, $H=(i\nabla+A^{(tr)}(x))^2$.
Set
\[
u(\omega)=U(\omega)-U(-\omega)
\]
where the function  $U(x)$
is defined by formula $(\ref{eq:Curv})$.
Denote by $S_0$   the operator of
multiplication by the function $\exp(iu(\omega) )$.
Then  all the results of
Theorem~$\ref{main}$ about the operator  $T (\lambda)$ are true for the operator
$T^{(tr)}(\lambda)= S^{(tr)}(\lambda)-S_0$. 
 \end{proposition}
 
 If a  magnetic field $B(x)$ satisfies assumption (\ref{eq:B1x}), then, similarly to the two-dimensional case (see subsection~4.1), Proposition~\ref{B6tr} can be deduced from Theorem~\ref{SM4} and Proposition~\ref{SpHomXY}. Such approach was used in
 \cite{R}. Indeed, if ${\rm curl}\, A(x)=o(|x|^{-2})$ as
$|x|\rightarrow\infty$, then necessarily condition (\ref{eq:reg1}) is satisfied
and hence function (\ref{eq:Curv}) is correctly defined. In this case
  the function $I(x,\omega)$ does not   depend on $x$
and equals $u(\omega)$.
Actually, it follows from (\ref{eq:GG}) that
\begin{eqnarray}\label{eq:GG2}
I(x,\omega)=\int_{-\infty}^\infty \langle{\rm
grad}\,U(x+t\omega),\omega\rangle d t 
= \lim_{T\rightarrow\infty}\int_{-T}^T \frac{d}{dt}U(x+t\omega) dt
 \nonumber\\ 
=\lim_{T\rightarrow\infty}(U(x+T\omega)-U(x-T\omega)) =U(\omega)-U(-\omega).
\end{eqnarray}
Therefore   function (\ref{eq:qqq1A}) equals zero
  so that the singular integral operator
  disappears in (\ref{eq:ssing}). 

 In the dimension $2$ the construction above works if (and only if) the total magnetic flux $\Phi$ is zero. Indeed, in this case
 \[
 \int_{|x|=R} \langle A^{(\infty)}(x), dx\rangle =0
 \]
 for any $R> 0$ so that   function (\ref{eq:Curv}) is  again correctly defined. Then Proposition~\ref{pot} for potential
 (\ref{eq:GGx}) and Theorem~\ref{main} for the SM remain true. The only difference is that under assumption (\ref{eq:B1x}) the integral    kernel of the   operator $T(\lambda)$ is $O(|\omega-\omega^\prime|^{-3+ r })$ 
as $\omega-\omega^\prime \rightarrow 0$.  
  
  \bigskip

{\bf 5.2.}
 As a concrete   example, let us   consider   a toroidal solenoid
 ${\bf T}$ in the space ${\Bbb R}^3$  symmetric with respect to rotations around the $x_3$-axis  
(which does not intersect  ${\bf T}$).  Suppose
 (which looks quite realistic) that a magnetic field
 \[
 B(x_1,x_2,x_3)= - \alpha (x_1^2+x_2^2)^{-1}(-x_2, x_1,0),  \quad
\alpha={\rm const},
\]
 inside of ${\bf T}$ and is zero outside. Then ${\rm div}\,B(x)=0$ and
the current ${\rm curl}\,B(x)=0$  if $x\not\in\partial{\bf T}$.
Of course, Theorem~\ref{main} applies to this field and hence
  $ S(\lambda)-I$   is integral operator with kernel from the class
$C^\infty({\Bbb S}^2\times {\Bbb S}^2)$. 

Let us illustrate our construction on this example.
First, we construct the potential $A^{(tr)}(x)$ by formula (\ref{eq:BB2}).
  We   assume that  the section ${\bf S}$ of ${\bf T}$, for example,  by  the half-plane   $x_2=0$, $x_1\geq 0$ is strictly convex
and has a smooth boundary
$\partial{\bf S}$ but is not necessarily  a disc.
Let the half-line $L_z$, $ z\in {\Bbb R}$,  consist of points $s (1+z^2)^{-1/2}(1,0,z)$ for all $s\in{\Bbb
R}_+$. Denote by  $z_1$ and  $z_2$ the values of $z$ for which   $L_z$ 
 is tangent to $\partial{\bf S}$ and, for $z\in[z_1,z_2]$,  
denote by $\varkappa_\pm(z)$, $\varkappa_+(z)\geq
\varkappa_-(z)$, the values of $s$ for which $L_z$ intersects $  {\bf S}$.
For $x=(x_1,x_2,x_3)$,  set $z=z(x)=x_3 (x_1^2+ x_2^2)^{-1/2}$. Taking into account the rotational symmetry, we
see that a point $s x \in{\bf T}$ if and only if
$s|x|(1+z^2)^{-1/2}(1,0,z)\in{\bf S}$ or $\varkappa_-(z)\leq s|x| \leq
\varkappa_+(z)$. Thus, integral (\ref{eq:BB2}) equals zero (and hence $A^{(tr)}(x)=0$)
if $z(x)\not\in (z_1,z_2)$ and it is actually taken over the set 
$[0,1]\cap [\varkappa_-(z)|x|^{-1},
\varkappa_+(z)|x|^{-1}]$ if  $z(x) \in (z_1,z_2)$.
Therefore $A^{(tr)}(x)=0$ if $|x|\leq \varkappa_-(z)$,
\begin{eqnarray}\label{eq:lagrA1}
 A^{(tr)}_1(x)&=& -\alpha x_1 x_3 (x_1^2+x_2^2)^{-1}\int_{\varkappa_-(z)/|x|}^1 ds
  \nonumber\\
   &=&-\alpha x_1 x_3
(x_1^2+x_2^2)^{-1}(1-\varkappa_-(z)/|x|)=:A^{(0)}_1(x) 
\end{eqnarray}
 if $\varkappa_-(z)\leq |x|\leq \varkappa_+(z)$ and
\begin{eqnarray}\label{eq:lagrA2}
 A^{(tr)}_1(x)&=& -\alpha x_1 x_3 (x_1^2+x_2^2)^{-1}\int_{\varkappa_-(z)/|x|}^{\varkappa_+(z)/|x|} ds
     \nonumber\\
 &=& -\alpha x_1 x_3 (x_1^2+x_2^2)^{-1}|x|^{-1}(\varkappa_+(z)-\varkappa_-(z) )  
\end{eqnarray}
 if $  |x|\geq \varkappa_+(z)$. Components $A^{(tr)}_2(x)$ and $A^{(tr)}_3(x)$ can be  found quite
similarly. In particular,
\begin{equation} \left.\begin{array}{lcl}
   A^{(\infty)}_1(x)&=& x_1 x_3 (x_1^2+x_2^2)^{-1} |x|^{-1} g
(z),
 \\
 A^{(\infty)}_2(x)&=& x_2 x_3 (x_1^2+x_2^2)^{-1} |x|^{-1} g (z),
 \\ 
A^{(\infty)}_3(x)&=&-  |x|^{-1} g (z), 
\end{array}\right\}
\label{eq:lagrA}\end{equation} 
  where
\[ 
 g (z)=- \alpha(\varkappa_+ (z) - \varkappa_- (z) ).
\]
 Clearly, $g(z)$   is  a continuous function, $\pm g(z)> 0$ if $\mp\alpha > 0$ for
$z \in(z_1,z_2)$ and $g(z)=0$   for $z\not\in (z_1,z_2)$. 
  
  Let ${\bf K}$ be the cone in ${\Bbb R}^3$ where $z (x) \in [z_{1}, z_{2}]$. Then
  ${\bf T}\subset {\bf K}$,  and ${\bf T}$ and $ {\bf K}$ are tangent to each other.
  The internal (external) part of ${\bf K}\setminus {\bf T}$
  will be denoted ${\bf K}_{int}$ (${\bf K}_{ext}$).  Of course $A^{(tr)}(x)=0$  if $x\not\in{\bf K} $. It follows from  (\ref{eq:lagrA1}), (\ref{eq:lagrA2}) that
    \[
    \left.\begin{array}{lcl}
   A^{(tr)}(x)&=&0, \quad\quad \quad\quad x\in{\bf K}_{int},
 \\
 A^{(tr)}(x)&=&A^{(0)}(x), \quad\quad x\in{\bf T},
 \\ 
A^{(tr)}(x)&=&A^{(\infty)}(x), \quad  x\in{\bf K}_{ext}. 
 \end{array}\right\}
    \] 
    Now formula (\ref{eq:B3}) for $A^{(tr)}(x)$ implies that
\begin{equation}\label{eq:Areg}
 \left.\begin{array}{lcl} 
A^{(reg)}(x)&=&-A^{( \infty)}(x) , \qquad\qquad  x\in{\bf K}_{int},
 \\
A^{(reg)}(x)&=&A^{(0)}(x)-A^{( \infty)}(x), \quad x\in{\bf T},
 \\ 
A^{(reg)}(x)&=&0, \qquad \qquad\qquad  x\in{\bf K}_{ext}. 
\end{array}\right\}
 \end{equation}

 Taking into account (\ref{eq:lagrA}),
 we see that a function
$U(x)$ satisfying (\ref{eq:GG}) can be constructed by the explicit formula
\begin{equation}\label{eq:UU1}
 U(x)=G(x_3 (x_1^2+ x_2^2)^{-1/2}),
\end{equation}
 where
\begin{equation}\label{eq:UU2}
 G^\prime (z)=-g(z)(z^2+1)^{-1/2}.
\end{equation}
 In particular, we see that $U(x)$ is a constant for $z (x)\not\in (z_{1},z_{2})$. Since
 $U(x)$ is defined up to a constant, we can set $U(x)=0$ for $z(x)\leq z_{1}$. Then
 \[
 U(x)=U_{0}=-\int_{- \infty}^\infty g(t)(t^2+1)^{-1/2} dt
\]
 for $z(x)\geq z_{2}$. It is easy to check that $-U_{0}$
  equals the
magnetic flux $\Phi_s$ through the section ${\bf S}$ of the   solenoid ${\bf T}$.
Indeed, let  $\omega_0=(0,0,1)$, $\langle x_0,\omega_0\rangle=0$,   and let
  $|x_0|$ and $R$ be sufficiently large.  By the Stokes  theorem,
$\Phi_s$ equals
the circulation of the   potential
$A^{(tr)}(x)$  over the closed contour formed by the four intervals
  $(-R\omega_{0}, R\omega_{0} )$, $(R\omega_{0}, R\omega_{0} +x_{0})$,
  $(  R\omega_{0} +x_{0}, -R\omega_{0}+x_{0})$ and $( - R\omega_{0} +x_{0}, -R\omega_{0})$.    Remark that
  $A^{(tr)}(x)\neq 0$ only on the interval $(  R\omega_{0} +x_{0}, -R\omega_{0}+x_{0})$ where $A^{(tr)}(x)=A^{(\infty)}(x)$, so that
  \[
 \Phi_{s} = - \int_{-R}^R \langle  A^{(\infty)}(x_{0}+t\omega_{0}), \omega_{0}\rangle dt.
 \]
   Passing to the limit $R\rightarrow\infty$ and using (\ref{eq:EABAB}), we see that $\Phi_s=-I(x_0,\omega_0)$. Hence equality (\ref{eq:GG2}) implies that
  \begin{equation}\label{eq:UU3}
\Phi_{s} =-U(\omega_{0})=-U_{0}.
\end{equation}

   Suppose now that the number $R_{2}$ in the definition of the cut-off function $\eta(x)$ is chosen in such a way that the ball $|x|\leq R_{2}$ does not intersect ${\bf T} $.
   Let first $x\not\in{\bf K} $ so that  $A^{(tr)}(x)=0$. Then the first formula   (\ref{eq:GGx}) shows  that $A(x)=0$  if   $z(x)\leq z_{1}$, and
\[
A(x) =- U_{0}\,  {\rm grad}\,  \eta(x), \quad x\not\in{\bf K}, \quad z(x) \geq z_{2}.
\] 
   If $x\in{\bf K} $, then according to the second formula   (\ref{eq:GGx}) and
   equalities (\ref{eq:Areg})
 \begin{eqnarray*}  
A(x)&=&0,  \qquad \qquad\qquad \qquad |x|\leq R_{1},
\nonumber\\
A(x)&=& -{\rm grad}\, ( \eta(x) U(x)),  \qquad R_{1}\leq |x|\leq R_{2},
\nonumber\\ 
A(x) &=& A^{(reg)}(x), \qquad\qquad\qquad |x|\geq R_{2}. 
 \end{eqnarray*}  
    In particular, $A(x)=0$ if $x\in{\bf T}_{ext}$.

The function $g(z)$ can be calculated explicitly if   ${\bf S}$     is  a disc. Suppose that this disc has radius
$r$, its center belongs to the    $x_3 $-axis and the distance from the center to  the $x_3$-axis
is $l$, $l>r$. Then the equation of   $\partial{\bf T}$   is
\begin{equation}\label{eq:af}
((x_1^2 +x_2^2)^{1/2}-l)^2+x_3^2=r^2.
\end{equation}
    Setting here $x_2=0$, $x_3=zx_1$, we obtain an equation
 for $x_1=x_1(z)$. The roots of this equation yield us
the numbers $(1+z^2)^{-1/2}\varkappa_\pm (z)$. Thus, 
\[
\varkappa_\pm (z)=(1+z^2)^{-1/2}(l\pm (r^2-(l^2-r^2)z^2)^{1/2})
\]
 and hence
\[
g(z)=-2\alpha(1+z^2)^{-1/2} (r^2-(l^2-r^2)z^2)^{1/2}.
\]
 In particular, $-z_1=z_2=r(l^2-r^2)^{-1/2}$ for this function.

Returning to the general case, we emphasize that a potential $A(x)$ satisfying the conclusions of Proposition~\ref{pot} is highly non-unique. Actually, the gradient of an arbitrary short-range function can be added to $A(x)$. For example, in the case (\ref{eq:af}) the magnetic potential completely different from the one constructed above can be found in the book \cite{Af}.
  
Let us finally calculate the SM $S^{(tr)}(\lambda)$.
  According to Proposition~\ref{B6tr},  up to an integral operator with
$C^\infty$-kernel, the SM $S^{(tr)}(\lambda)$ is  the operator
$S_0$ of multiplication by the function $\exp (iu(\omega))$,
where by virtue of (\ref{eq:UU1}), (\ref{eq:UU2})
\[
u(\omega)= q(\omega_3 (1- \omega_3^2)^{-1/2}) ,\quad  
\omega =(\omega_1,\omega_2,\omega_3)\in {\Bbb S}^2,
\] 
and
\[
 q(z)=G(z)-G(-z)=-\int_{-z}^z g(t)(t^2+1)^{-1/2} dt.
\]
 Clearly, $q(-z)=-q(z) $, $ q( z) $ is an increasing (decreasing)
function if $\alpha > 0$ ($\alpha< 0$) and it
is a constant, $ q( z)=q_0 $, if
$z\geq \max\{ |z_1|, |z_2| \}$. It follows from (\ref{eq:UU3}) that
\[
q_0= u(\omega_0)= U (\omega_0)- U (-\omega_0)=
-\Phi_s.
\]  
  Thus, $u(\omega)$ depends only on the coordinate $\omega_{3}$ and takes all the values between $-|\Phi_s|$ and $|\Phi_s|$. Therefore $\sigma_{ess}(S^{(tr)}(\lambda))$ 
   coincides with the arc $[e^{-i|\Phi_s|}, e^{i|\Phi_s|}]$ if $|\Phi_s| < \pi$,  and it  covers the whole unit circle if $|\Phi_s|\geq \pi$.

\end{document}